# GENERALIZING BENFORD'S LAW USING POWER LAWS:

# APPLICATION TO INTEGER SEQUENCES

Hürlimann Werner, Feldstrasse 145
CH-8004 Zürich, Switzerland
E-mail : whurlimann@bluewin.ch
URL : www.geocities.com/hurlimann53

**Abstract.**

A simple method to derive parametric analytical extensions of Benford's law for first digits of numerical data is proposed. Two generalized Benford distributions are considered, namely the two-sided power Benford (TSPB) distribution, which has been introduced in Hürlimann(2003), and the new Pareto Benford (PB) distribution. Based on the minimum chi-square estimators, the fitting capabilities of these generalized Benford distributions are illustrated and compared at some interesting and important integer sequences. In particular, it is significant that much of the analyzed integer sequences follow with a high p-value the generalized Benford distributions. While the sequences of prime numbers less than 1,000 respectively 10,000 are not at all Benford or TSPB distributed, they are approximately PB distributed with high p-values of 93.3% and 99.9%. Benford's law of a mixing of data sets is rejected at the 5% significance level while the PB law is accepted with a 93.6% p-value, which improves the p-value of 25.2%, which has been obtained previously for the TSPB law.

## 1. Introduction.

Since Newcomb(1881) and Benford(1938) it is known that many numerical data sets follow Benford's law or are closely approximated by it. To be specific, if the random variable $X$, which describes the first significant digit in a numerical table, is Benford distributed, then

$$P(X = d) = \log(1 + d^{-1}), \quad d \in \{1,...,9\}. \tag{1.1}$$

Mathematical explanations of this law have been proposed by Pinkham(1961), Cohen(1976), Hill(1995a/b/c,97,98), Allart(1997), Janvresse and de la Rue(2004). In recent years an upsurge of applications of Benford's law has appeared, as can be seen from the recently compiled bibliography by Hürlimann(2006).

Hill(1995c) also suggested to switch the attention to probability distributions that follow or closely approximate Benford's law. Papers along this path include Leemis et al.(2000) and Engel and Leuenberger(2003). Some survival distributions, which satisfy exactly Benford's law, are known. However, many simple analytical distributions, which include as special case Benford's law are not known. Combining facts from Leemis et al.(2000) and Dorp and Kotz(2002), such a simple one-parameter family of distributions has been considered in Hürlimann(2003). In a sequel to this, a further generalization of Benford's law is considered.

The interest of enlarged Benford laws is two-fold. First, parametric extensions may provide a better fit of the data than Benford's law itself. Second, they yield a simple statistical procedure to validate Benford's law. If Benford's model is sufficiently "close" to the one-parameter extended model, then it will be retained. These points will be illustrated through our application to integer sequences.

## 2. <u>Generalizing Benford's distribution.</u>

If $T$ denotes a random lifetime with survival distribution $S(t) = P(T \geq t)$, then the value $Y$ of the first significant digit in the lifetime $T$ has probability distribution

$$P(Y = y) = \sum_{i=-\infty}^{\infty} \left\{ S(y \cdot 10^i) - S((y+1) \cdot 10^i) \right\}, \quad y \in \{1, \ldots, 9\}. \qquad (2.1)$$

Alternatively, if $D$ denotes the integer-valued random variable satisfying

$$10^D \leq T < 10^{D+1}, \qquad (2.2)$$

then the first significant digit can be written in terms of $T$ and $D$ as

$$Y = \left[ T \cdot 10^{-D} \right] = \left[ 10^{\log T - D} \right], \qquad (2.3)$$

where $[x]$ denotes the greatest integer less than or equal to $x$. In particular, if the random variable $Z = \log T - D$ is uniformly distributed as $U(0,1)$, then the first significant digit $Y$ is exactly Benford distributed. Starting from the *uniform* random variable $W = U(0,2)$ or the *triangular* random variable $W = Triangular(0,1,2)$ with probability density function $f_W(w) = w$ if $w \in (0,1)$ and $f_W(w) = 2 - w$ if $w \in [1,2)$, one shows that the random lifetime $T = 10^W$ generates the first digit Benford distribution (Leemis et al.(2000), Examples 1 and 2).

A simple parametric distribution, which includes as special cases both the above uniform and triangular distributions, is the two-sided power random variable $W = TSP(\alpha, c)$ considered in Dorp and Kotz(2002) with probability density function

$$f_W(w) = \begin{cases} \frac{c}{2} \left( \frac{w}{\alpha} \right)^{c-1}, & 0 < w \leq \alpha, \\ \frac{c}{2} \left( \frac{2-w}{2-\alpha} \right)^{c-1}, & \alpha \leq w < 2. \end{cases} \qquad (2.4)$$

If $c = 1$ then $W = U(0,2)$, and if $c = 2, \alpha = 1$ then $W = Triangular(0,1,2)$. This observation shows that the random lifetime $T = 10^{TSP(1,c)}$ will generate first digit distributions closely related to Benford's distribution, at least if $c$ is close to 1 or 2.

**<u>Theorem 2.1.</u>** Let $W = TSP(1,c)$ be the two-sided power random variable with probability density function

$$f_W(w) = \begin{cases} \frac{c}{2}w^{c-1}, & 0 < w \le 1, \\ \frac{c}{2}(2-w)^{c-1}, & 1 \le w < 2, \end{cases} \tag{2.5}$$

and let the integer-valued random variable $D$ satisfy $D \le W < D+1$. Then the first digit random variable $Y = \left[10^{W-D}\right]$ has the one-parameter *two-sided power Benford* (TSPB) probability density function

$$f_Y(y) = \frac{1}{2}\left\{\left[\log(1+y)\right]^c - \left[\log y\right]^c - \left[1-\log(1+y)\right]^c + \left[1-\log y\right]^c\right\}, \quad y \in \{1,...,9\}. \tag{2.6}$$

**Proof.** This has been shown in Hürlimann(2003). ◊

## 3.  **From the geometric Brownian motion to the Pareto Benford law.**

Another interesting distribution, which also takes the form of a two-sided power law, is the double Pareto random variable $W = DP(s,\alpha,\beta)$ considered in Reed(2001) with probability density function

$$f_W(w) = \begin{cases} \frac{\alpha\beta}{\alpha+\beta}\left(\frac{w}{s}\right)^{\beta-1}, & w \le s, \\ \frac{\alpha\beta}{\alpha+\beta}\left(\frac{w}{s}\right)^{-\alpha-1}, & w \ge s. \end{cases} \tag{3.1}$$

Recall the stochastic mechanism and the natural motivation, which generates this distribution. It is often assumed that the time evolution of a stochastic phenomena $X_t$ involves a variable but size independent proportional growth rate and can thus be modeled by a geometric Brownian motion (GBM) described by the stochastic differential equation

$$dX = \mu \cdot X \cdot dt + \sigma \cdot X \cdot dW, \tag{3.2}$$

where $dW$ is the increment of a Wiener process. Since the proportional increment of a GBM in time $dt$ has a systematic component $\mu \cdot dt$ and a random white noise component $\sigma \cdot dW$, GBM can be viewed as a stochastic version of a simple exponential growth model. The GBM has long been used to model the evolution of stock prices (Black-Scholes option pricing model), firm sizes, city sizes and individual incomes. It is well-known that empirical studies on such phenomena often exhibit power-law behavior. However, the state of a GBM after a fixed time $T$ follows a lognormal distribution, which does not exhibit power-law behavior.

Why does one observe power-law behavior for phenomena apparently evolving like a GBM? A simple mechanism, which generates the power-law behavior in the tails, consists to assume that the time of observation $T$ itself is a random variable, whose distribution is an exponential distribution. The distribution of $X_T$ with fixed initial state $s$ is described by the *double Pareto distribution* $DP(s,\alpha,\beta)$ with density function (3.1), where $\alpha,\beta > 0$, and $\alpha, -\beta$ are the positive roots of the characteristic equation

$$\left(\mu - \frac{1}{2}\sigma^2\right)z + \frac{1}{2}\sigma^2 z^2 = \lambda, \tag{3.3}$$

where $\lambda$ is the parameter of the exponentially distributed random variable $T$. Setting $s = 1$ one obtains the following generalized Benford distribution.

**Theorem 3.1.** Let $W = DP(1, \alpha, \beta)$ be the double Pareto random variable with probability density function

$$f_W(w) = \begin{cases} \frac{\alpha\beta}{\alpha+\beta}(w)^{\beta-1}, & w \leq 1, \\ \frac{\alpha\beta}{\alpha+\beta}(w)^{-\alpha-1}, & w \geq 1. \end{cases} \qquad (3.4)$$

Let the integer-valued random variable $D$ satisfy $D \leq W < D+1$. Then the first digit random variable $Y = \left[10^{W-D}\right]$ has the two-parameter *Pareto Benford* (PB) probability density function

$$f_Y(y) = \frac{\alpha}{\alpha+\beta}\left\{[\log(1+y)]^\beta - [\log(y)]^\beta\right\}$$
$$+ \frac{\beta}{\alpha+\beta}\cdot\sum_{k=1}^{\infty}\left\{[k+\log(y)]^{-\alpha} - [k+\log(1+y)]^{-\alpha}\right\}, \quad y \in \{1,...,9\}. \qquad (3.5)$$

**Proof.** The probability density function of $T = 10^W$ is given by

$$f_T(t) = \frac{1}{t \cdot \ln 10}\cdot f_W\left(\frac{\ln t}{\ln 10}\right) = \begin{cases} \frac{\alpha\beta}{\alpha+\beta}\frac{1}{t\cdot\ln 10}\left(\frac{\ln t}{\ln 10}\right)^{\beta-1}, & 1 < t \leq 10, \\ \frac{\alpha\beta}{\alpha+\beta}\frac{1}{t\cdot\ln 10}\left(\frac{\ln t}{\ln 10}\right)^{-\alpha-1}, & t > 10. \end{cases}$$

It follows that the first significant digit of $T$, namely $Y = \left[T \cdot 10^{-D}\right]$, has probability density

$$f_Y(y) = \sum_{k=0}^{\infty}\int_{10^k y}^{10^k(y+1)} f_T(t)dt.$$

Making the change of variable $u = \ln t / \ln 10$, one obtains

$$f_Y(y) = \frac{\alpha\beta}{\alpha+\beta}\left\{\int_{\log y}^{\log(y+1)} u^{\beta-1}du + \sum_{k=1}^{\infty}\int_{k+\log y}^{k+\log(y+1)} u^{-\alpha-1}du\right\} = \frac{\alpha\beta}{\alpha+\beta}\left\{\frac{1}{\beta}u^\beta\left|\begin{matrix}\log(1+y)\\\log(y)\end{matrix}\right. + \sum_{k=1}^{\infty}\frac{-1}{\alpha}u^{-\alpha}\left|\begin{matrix}k+\log(1+y)\\k+\log(y)\end{matrix}\right.\right\}$$

$$= \frac{\alpha}{\alpha+\beta}\left\{[\log(1+y)]^\beta - [\log(y)]^\beta\right\} + \frac{\beta}{\alpha+\beta}\cdot\sum_{k=1}^{\infty}\left\{[k+\log(y)]^{-\alpha} - [k+\log(1+y)]^{-\alpha}\right\},$$

which is (3.5). ◊

One observes that setting $\beta = 1$ and letting $\alpha$ goes to infinity, the Pareto Benford distribution converges to Benford's law.

## 4. Fitting the first digit distributions of integer sequences.

Minimum chi-square estimation of the generalized Benford distributions is straightforward by calculation with modern computer algebra systems. The fitting capabilities of the new distributions are illustrated at some interesting and important integer sequences.

The first digit occurrences of the analyzed integer sequences are listed in Table 4.1. The minimum chi-square estimators of the generalized distributions as well as an assumed summation index $m$ for the infinite series (3.5) are displayed in Table 4.2. Statistical results are summarized in Table 4.3. For comparison we list the chi-square values and their corresponding p-values. The obtained results are discussed.

The definition, origin and comments on the mathematical interest of a great part of these integer sequences has been discussed in Hürlimann(2003). Further details on all sequences can be retrieved from the considerable related literature. The *mixing sequence* represents the aggregate of the integer sequences considered in Hürlimann(2003).

**Table 4.1:** First digit distributions of some integer sequences

| Name of sequence | Sample size | Percentage of first digit occurrences | | | | | | | | |
|---|---|---|---|---|---|---|---|---|---|---|
| | | 1 | 2 | 3 | 4 | 5 | 6 | 7 | 8 | 9 |
| Benford law | | 30.1 | 17.6 | 12.5 | 9.7 | 7.9 | 6.7 | 5.8 | 5.1 | 4.6 |
| Square | 100 | 21.0 | 14.0 | 12.0 | 12.0 | 9.0 | 9.0 | 8.0 | 7.0 | 8.0 |
| Cube | 500 | 28.2 | 14.8 | 11.4 | 9.8 | 8.8 | 7.8 | 6.6 | 6.8 | 5.8 |
| Cube | 1000 | 22.6 | 15.9 | 12.4 | 10.6 | 9.4 | 8.3 | 7.4 | 7.1 | 6.3 |
| Cube | 10000 | 22.5 | 15.8 | 12.6 | 10.6 | 9.3 | 8.3 | 7.5 | 7.0 | 6.4 |
| Square root | 99 | 19.2 | 17.2 | 15.2 | 13.1 | 11.1 | 9.1 | 7.1 | 5.1 | 3.0 |
| Prime < 100 | 25 | 16.0 | 12.0 | 12.0 | 12.0 | 12.0 | 8.0 | 16.0 | 8.0 | 4.0 |
| Prime < 1000 | 168 | 14.9 | 11.3 | 11.3 | 11.9 | 10.1 | 10.7 | 10.7 | 10.1 | 8.9 |
| Prime < 10000 | 1229 | 13.0 | 11.3 | 11.9 | 11.3 | 10.7 | 11.0 | 10.2 | 10.3 | 10.3 |
| Princeton number | 25 | 28.0 | 8.0 | 12.0 | 12.0 | 8.0 | 12.0 | 8.0 | 4.0 | 8.0 |
| Mixing sequence | 618 | 28.3 | 14.6 | 11.5 | 9.9 | 7.6 | 7.8 | 8.1 | 6.6 | 5.7 |
| Pentagonal number | 100 | 35.0 | 12.0 | 10.0 | 8.0 | 10.0 | 6.0 | 8.0 | 5.0 | 6.0 |
| Keith number | 71 | 32.4 | 14.1 | 14.1 | 7.0 | 4.2 | 7.0 | 12.7 | 2.8 | 5.6 |
| Bell number | 100 | 31.0 | 15.0 | 10.0 | 12.0 | 10.0 | 8.0 | 5.0 | 6.0 | 3.0 |
| Catalan number | 100 | 33.0 | 18.0 | 11.0 | 11.0 | 8.0 | 8.0 | 4.0 | 3.0 | 4.0 |
| Lucky number | 45 | 42.2 | 17.8 | 8.9 | 4.4 | 2.2 | 6.7 | 8.9 | 2.2 | 6.7 |
| Ulam number | 44 | 45.5 | 13.6 | 6.8 | 6.8 | 4.5 | 6.8 | 4.5 | 6.8 | 4.5 |
| Numeri ideoni | 65 | 30.8 | 18.5 | 13.8 | 10.8 | 6.2 | 3.1 | 7.7 | 6.2 | 3.1 |
| Fibonacci number | 100 | 30.0 | 18.0 | 13.0 | 9.0 | 8.0 | 6.0 | 5.0 | 7.0 | 4.0 |
| Partition number | 94 | 28.7 | 17.0 | 14.9 | 9.6 | 7.4 | 6.4 | 7.4 | 5.3 | 3.2 |

**Table 4.2:** Minimum chi-square estimators

| Name of sequence | Sample size | TSPB parameter | PB parameters | | |
|---|---|---|---|---|---|
| | | c | alpha | beta | m |
| Square | 100 | 0.79837 | 15.55957 | 1.74552 | 100 |
| Cube | 500 | 2.46519 | 5.55849 | 1.69860 | 100 |
| Cube | 1000 | 2.26798 | 20.56506 | 1.47082 | 100 |
| Cube | 10000 | 2.27054 | 20.53577 | 1.47576 | 100 |
| Square root | 99 | 1.40176 | 89491723 | 1.34334 | 100 |
| Prime < 100 | 25 | 2.68581 | 23.13952 | 2.14449 | 100 |
| Prime < 1000 | 168 | 2.95216 | 22.99754 | 2.28436 | 100 |
| Prime < 10000 | 1229 | 3.03542 | 29.76729 | 2.30760 | 100 |
| Princeton number | 25 | 2.76170 | 6.94595 | 2.36119 | 100 |
| Mixing sequence | 618 | 2.53958 | 4.78641 | 1.83119 | 100 |
| Pentagonal number | 100 | 2.94847 | 2.06797 | 3.31268 | 100 |
| Keith number | 71 | 2.73338 | 2.16107 | 2.63720 | 1000 |
| Bell number | 100 | 1.08191 | 10.14820 | 1.24828 | 100 |
| Catalan number | 100 | 1.13522 | 0.67095 | 1.15377 | 5000 |
| Lucky number | 45 | 3.15721 | 7.56962 | 0.94576 | 100 |
| Ulam number | 44 | 3.55375 | 9.99445 | 0.81215 | 100 |
| Numeri ideoni | 65 | 1.12410 | 1297612.16 | 0.98591 | 100 |
| Fibonacci number | 100 | 2.05365 | 257000.42 | 1.00560 | 100 |
| Partition number | 94 | 1.23268 | 0.65651 | 1.71409 | 1000 |

**Table 4.3:** Fitting integer sequences to the Benford and generalized Benford distributions

| Name of sequence | Sample size | Benford chi-square | P-value | Two-Sided Power Benford chi-square | P-value | Pareto Benford chi-square | P-value |
|---|---|---|---|---|---|---|---|
| Square | 100 | 9.096 | 33.43 | 7.837 | 34.72 | 0.362 | 99.91 |
| Cube | 500 | 9.696 | 28.70 | 5.808 | 56.23 | 0.286 | 99.96 |
| Cube | 1000 | 46.459 | 0.00 | 43.725 | 0.00 | 0.480 | 99.81 |
| Cube | 10000 | 443.745 | 0.00 | 472.011 | 0.00 | 3.138 | 79.13 |
| Square root | 99 | 8.612 | 37.61 | 7.002 | 42.86 | 2.778 | 83.61 |
| Prime < 100 | 25 | 7.741 | 45.91 | 7.299 | 39.84 | 1.849 | 93.30 |
| Prime < 1000 | 168 | 45.016 | 0.00 | 36.651 | 0.00 | 0.333 | 99.93 |
| Prime < 10000 | 1229 | 387.194 | 0.00 | 307.322 | 0.00 | 3.297 | 77.07 |
| Princeton number | 25 | 3.452 | 90.29 | 2.762 | 89.72 | 1.302 | 97.16 |
| Mixing sequence | 618 | 15.550 | 4.93 | 9.014 | 25.17 | 1.819 | 93.55 |
| Pentagonal number | 100 | 5.277 | 72.76 | 2.127 | 95.24 | 1.968 | 92.26 |
| Keith number | 71 | 9.215 | 32.45 | 7.688 | 36.09 | 7.402 | 28.53 |
| Bell number | 100 | 3.069 | 93.00 | 3.014 | 88.37 | 2.607 | 85.63 |
| Catalan number | 100 | 2.404 | 96.61 | 2.304 | 94.11 | 1.934 | 92.57 |
| Lucky number | 45 | 7.693 | 46.40 | 5.165 | 63.98 | 5.564 | 47.37 |
| Ulam number | 44 | 6.350 | 60.81 | 2.520 | 92.56 | 2.526 | 86.56 |
| Numeri ideoni | 65 | 2.594 | 95.72 | 2.522 | 92.54 | 2.584 | 85.89 |
| Fibonacci number | 100 | 1.029 | 99.81 | 1.021 | 99.45 | 1.027 | 98.46 |
| Partition number | 94 | 1.394 | 99.43 | 1.132 | 99.24 | 1.513 | 95.86 |

All of the 19 considered integer sequences are quite well fitted by the new PB distribution. For 14 sequences the minimum chi-square is smallest among the three comparative values and in the other 5 cases its value does not differ much from the chi-square of the TSPB distribution (green cells in Table 4.3).

The Benford property for the sequence of primes has long been studied. Diaconis(1977) shows that primes are not Benford distributed. However, it is known that the sequence of primes is Benford distributed with respect to other densities rather than with the usual natural density (Whitney(1972), Schatte(1983), Cohen and Katz(1983)). Bombieri (see Serre(1996), p. 76) has noted that the analytical density of primes with first digit 1 is $\log_{10} 2$, and this result can be easily generalized to Benford behavior for any first digit. Table 4.3 shows that the primes less than 1,000 respectively 10,000 are not at all Benford or TSPB distributed, but they are approximately PB distributed with high p-values of 93.3% and 99.9%. Is this a new property of the prime number sequence? Unfortunately, the fit of the Pareto Benford distribution for the 78,499 prime numbers below 100,000 is rejected since the corresponding minimum chi-square value equals 1,391. Therefore it seems that the good fit of the BP distribution remains restricted to finite prime number sequences. Similar results for the sequence of squares and cubes can be made. Recall that the exact probability distribution of the first digit of $m$-th integer powers with at most $n$ digits is known and asymptotically related to Benford's law (e.g. Hürlimann(2004)). Here again the fit of the PB distribution is very good when restricted to finite sequences but breaks down for longer sequences. A further remarkable result is that Benford's law of the mixing sequence is rejected at the 5% significance level while the PB law is accepted with a 93.6% p-value, which improves the p-value of 25.2% obtained previously for the TSPB law in Hürlimann(2003).

A strong numerical evidence for the Benford property for the Fibonacci, Bell, Catalan and partition numbers is observed (corresponding yellow cells in Tables 4.2 and 4.3). In particular, the values of the parameters $\alpha, \beta$ of the BP distribution for the Fibonacci sequence are close to 1 and $\infty$, which means that the BP distribution is almost Benford as remarked after Theorem 3.1. It is well-known that the Fibonacci sequence is Benford distributed (e.g. Brown and Duncan(1970), Wlodarski(1971), Sentence(1973), Webb(1975), Raimi(1976), Brady(1978) and Kunoff(1987)).The same result for Bell numbers has been derived formally in Hürlimann(2003), Theorem 4.1.

## References.

*Allaart, P.C.* (1997). An invariant-sum characterization of Benford's law. Journal of Applied Probability 34, 288-291.

*Benford, F.* (1938). The law of anomalous numbers. Proceedings of the American Philosophical Society 78, 551-572.

*Brady, W.G.* (1978). More on Benford's law. Fibonacci Quarterly 16, 51-52.

*Brown, J. and R. Duncan* (1970). Modulo one uniform distribution of the sequence of logarithms of certain recursive sequences. Fibonacci Quarterly 8, 482-486.

*Cohen, D.* (1976). An explanation of the first digit phenomenon. Journal of Combinatorial Theory (A) 20, 367-370.

*Cohen, D. and T. Katz* (1984). Prime numbers and the first digit phenomenon. Journal of Number Theory 18, 261-268.

*Dorp, J.R. van and S. Kotz* (2002). The standard two-sided power distribution and its properties : with applications in financial engineering. The American Statistician 56(2), 90-99.

*Engel, H.-A. and C. Leuenberger* (2003). Benford's law for exponential random variables. Statistics and Probability Letters 63(4), 361-365.

*Hill, T.P.* (1995a). Base-invariance implies Benford's law. Proceedings of the American Mathematical Society 123, 887-895.

*Hill, T.P.* (1995b). The significant-digit phenomenon. Amer. Math. Monthly 102, 322-326.

*Hill, T.P.* (1995c). A statistical derivation of the significant-digit law. Statistical Science 10, 354-363.

*Hill, T.P.* (1997). Benford's law. Encyclopedia of Mathematics Supplement, vol. 1, 102.

*Hill, T.P.* (1998). The first digit phenomenon. The American Scientist 10(4), 354-363.

*Hürlimann, W.* (2003). A generalized Benford law and its application. Advances and Applications in Statistics 3(3), 217-228.

*Hürlimann, W.* (2004). Integer powers and Benford's law. International Journal of Pure and Applied Mathematics 11(1), 39-46.

*Hürlimann, W.* (2006). Benford's law from 1881 to 2006: a bibliography. http://arxiv.org.

*Janvresse, E. and T. de la Rue* (2004). From uniform distribution to Benford's law. Journal of Applied Probability 41(4), 1203-1210.

*Kunoff, S.* (1987). N! has the first digit property. Fibonacci Quarterly 25, 365-367.

*Leemis, L.M., Schmeiser, B.W. and D.L. Evans* (2000). Survival distributions satisfying Benford's law. The American Statistician 54(3), 1-6.

*Newcomb, S.* (1881). Note on the frequency of use of the different digits in natural numbers. American Journal of Mathematics 4, 39-40.

*Pinkham, R.S.* (1961). On the distribution of first significant digits. Annals of Mathematical Statistics 32, 1223-1230.

*Raimi, R.A.* (1969). The peculiar distribution of first digits. Scientific American, 109-120.

*Reed, W. J.* (2001). The Pareto, Zipf and other power laws. Economics Letters 74, 15-19.

*Schatte, P.* (1983). On $H_\infty$-summability and the uniform distribution of sequences. Math. Nachr. 113, 237-243.

*Sentance, W.A.* (1973). A further analysis of Benford's law. Fibonacci Quarterly 11, 490-494.

*Serre, J.-P.* (1996). A Course in Arithmetic. Springer-Verlag.

*Webb, W.* (1975). Distribution of the first digits of Fibonacci numbers. Fibonacci Quarterly 13, 334-336.

*Whitney, R.E.* (1972). Initial digits for the sequence of primes. American Mathematical Monthly 79, 150-152.

*Wlodarski, J.* (1971). Fibonacci and Lucas Numbers tend to obey Benford's law. Fibonacci Quarterly 9, 87-88.